\documentclass[11pt,english,twoside,a4paper]{article}
\usepackage{latexsym}
\usepackage{babel}
\usepackage{amsmath}
\usepackage{amsfonts}
\usepackage{amssymb}
\usepackage{graphicx}
\usepackage{xcolor}
\DeclareFontFamily{OT1}{pzc}{}
\DeclareFontShape{OT1}{pzc}{m}{it}{<-> s * [1.200] pzcmi7t}{}
\DeclareMathAlphabet{\mathpzc}{OT1}{pzc}{m}{it}
\newcommand{\enne}{\mathbb{N}}
\newcommand{\erre}{\mathbb{R}}
\newcommand{\bv}{\mbox{\bf v}}
\newcommand{\bw}{\mbox{\bf w}}
\newcommand{\bg}{\bigskip}
\newcommand{\md}{\medskip}
\newcommand{\sm}{\smallskip}
\newcommand{\np}{\noindent}
\newcommand{\beq}{\begin{equation}}
\newcommand{\eeq}{\end{equation}}

\newcommand{\y}{\mathpzc{y}}
\newcommand{\red}{\textcolor{red}}

\newtheorem{theorem}{Theorem}[section]
\newtheorem{corollary}{Corollary}[section]
\newtheorem{proposition}{Proposition}[section]
\newtheorem{definition}{Definition}[section]
\newtheorem{remark}{Remark}[section]

\newtheorem{example}{Example}[section]
\title{
Uniform asymptotic stability of a PDE's system \\
arising from a flexible robotics model
}
\author{{\bf Tiziana Cardinali$^{(1)}$, Serena Matucci$^{(2)}$, Paola Rubbioni$^{(1)}$}\\
{\small {\em (1) Department of Mathematics and Computer Science, University of Perugia}}\\
{\small {\em (2) Department of Mathematics and Computer Sciences ``U. Dini'', University of Florence}}\\
{\small {\em E-mail addresses: 
tiziana.cardinali@unipg.it; serena.matucci@unifi.it; paola.rubbioni@unipg.it}}
}
\date{}
\begin{document}

\maketitle

\bg

\np
{\bf Abstract:} 
In this paper we investigate the asymptotic stability of a fourth-order PDE with a fading memory forcing term and boundary conditions arising from a flexible robotics model. We carry on our study by using an abstract formulation of the problem based on the $C_0$-semigroup. To achieve our objective, we first provide new results on the existence, uniqueness, continuous dependence on initial data of either mild and strong solutions for semilinear integro-differential equations in Banach spaces.  Then, we also find sufficient conditions for the uniform asymptotic stability of solutions and for the existence of attactors. As an application of these abstract results, we can ensure existence, uniqueness and continuous dependence on initial data for the solutions of the boundary value problem under investigation and, finally, we prove the uniform asymptotic stability of solutions and the existence of attactors under suitable conditions on the nonlinear term.
\bg

\np
{\em 2020 MSC:} Primary: 35G30;  34G20. 
Secondary: 35A24; 34D20. 
\md

\np
{\em Keywords and phrases.}  Buondary value problems; integro-differential equation; uniform asymptotic stability; mild solution; strong solution; flexible robotics.


\bg
\bg

\section{Introduction}
\label{s:i}

The present work is devoted to the study of the uniform asymptotic stability of solutions for a class of boundary value problems with controls arising from mechanical problems, and in particular from flexible beams and robotic models. 
The study of mechanical systems with flexible components is a current and interesting research topic, as several examples are found in robotics, soft robotics, materials science, and stretchable electronics.
The vibrations of a flexible beam with negligible inertia, clamped at one end and controlled at the free end, can be described by the Euler-Bernoulli beam equation
\beq
\label{e:pde0}\tag{E-B}
u_{tt}(t,x)+ u_{xxxx}(t,x)=0 ,
\eeq
with boundary conditions
\begin{eqnarray*}
&& u(t,0)=u_x(t,0)=u_{xx}(t,1)=0,\\
&& mu_{tt}(t,1)-u_{xxx}(t,1)= w(t), 
\end{eqnarray*}
for $t\in (0,+\infty)$, $x\in (0,1)$, where $u : \erre^+_0 \times [0,1] \to \erre$ is the vibration amplitude, $m$ is the tip mass, and $w$ is the boundary control force applied at the free end of the beam, often used to stabilize the system. 
\\
In 1998 Conrad and Morgül \cite{CM} proposed a kind of control given by the next linear feedback law 
\begin{equation*}
w(t)=-\alpha u_t(t,1) + \beta u_{xxxt}(t,1), \ t\in (0,+\infty),\ (\alpha, \beta >0)
\end{equation*}
and proved that the closed-loop system is well-posed and the solutions uniformly decay to zero.
\\
As far as we know, the Conrad and Morgül's  paper  is a first example of approaching by the semigroup theory the models described above. This method consists in transforming the model equation and boundary conditions into the linear ordinary differential equation of the first order
$$
y'(t)=Ay(t),\ \ t\in (0,+\infty),
$$
 in a suitable function space, with $A$ generator of a $C_0$-semigroup of contractions. 
\\
This approach has led to new results in several directions, see for instance the very recent works 
\cite{AOO,CF,DH,SBG}, and in particular in the field of stabilization of infinite-dimensional dynamic systems in abstract spaces.
 For these topics we also refer for example to the books \cite{LGM,SZ}.
\\
Further, in the last years several papers appeared where equation \eqref{e:pde0} is perturbed by different kinds of forcing terms. For instance, in \cite{GG}
the authors consider the equation
$$
u_{tt}(t,x) + u_{xxxx}(t,x)= - f(u)+g\left(t,x\right),
$$ 
where $f$ and $g$ are nonlinear forces, and study the instability of the system with boundary conditions to describe a model of suspended bridge. 
\sm

In this paper we deal with a model driven by the next fourth-order PDE including  a forcing term depending on a distributed delay, which formalizes the movement of a robotic arm having fading memory of the past deflections,
\begin{eqnarray}
\label{e:pde1}
&& u_{tt}(t,x)= -\, u_{xxxx}(t,x) + g\left(t,u(t,x),\int_{0}^t \frac{e^{-(t-s)/T}}{T}\, u(s,x)\, ds\right),
\end{eqnarray}
where  $u : \erre^+_0 \times [0,1] \to \erre$, $g : \erre^+_0 \times \erre \times \erre \to \erre$, $m,\alpha,\beta\in \erre^+$. 
The kernel
$$
k(t,s)=\frac{e^{-(t-s)/T}}{T}
$$
is given by the exponential distribution of probability ${\cal K}(r)=\dfrac{e^{-r/T}}{T}$, $r\ge 0$, describing a fading delay which can be regulated by the width $T$. 
\\
Equation \eqref{e:pde1} is here studied when subject to the boundary conditions
\begin{eqnarray}
\label{e:pde2}
&& u(t,0)=u_x(t,0)=u_{xx}(t,1)=0,\\
\label{e:pde3}
&& mu_{tt}(t,1)-u_{xxx}(t,1)=-\alpha u_t(t,1)+\beta u_{xxxt}(t,1),
\end{eqnarray}
where equation \eqref{e:pde3} represents a boundary control input at the free end of the arm, 
and to the initial data
\begin{eqnarray}
\label{e:iv0} && u(0,x)=\bar p(x),  \\
\label{e:iv1}&& u_t(0,x)=\bar q(x),\\
\label{e:iv2} && -u_{xxx}(0,1)+\frac{m}{\beta}u_t(0,1)=-\bar p_{xxx}(1)+\frac{m}{\beta}\bar q(1)=:\bar \eta,
\end{eqnarray}
for $t\in(0,+\infty), \, x\in (0,1)$, with $\bar p\in H^4(0,1)$ and $\bar q\in H^1(0,1)$ given functions. 
\\
As in \cite{CM}, without loss of generality, the lenght of the beam or arm, its  flexural rigidity, and the mass per unit lenght are chosen to be unity. 
\sm

Along the lines of \cite{CM}, we rewrite the problem \eqref{e:pde1}, \eqref{e:pde2}-\eqref{e:pde3}, \eqref{e:iv0}-\eqref{e:iv2} as a Cauchy problem in a suitable Hilbert space, driven by a semilinear integro-differential equation
(see Section \ref{s:abs}). 
\\
This equation can formally be seen of the type
\begin{equation}\label{e:eq}
y'(t)= Ay(t)+f\left(t,y(t),\int_{t_0}^t k(t,s)y(s) ds\right),\ t\ge t_0,
\end{equation}
and it can be studied in the more abstract setting of a Banach space $E$. 

We devote Section \ref{s:eudc} to establishing sufficient conditions for the existence of solutions, 
both in the mild and in the strong sense, to the initial value problem obtained coupling equation \eqref{e:eq} with $y(t_0)=\bv\in E$.
Theorems on the uniqueness and the continuous dependence on the initial data are provided as well.

These first theorems play a basic role for the results of the next section, which is the heart of the article. In Section \ref{s:as}, in fact, we provide one of the main theorems of this manuscript. We prove indeed the uniform asymptotical stability of the solutions of \eqref{e:eq}, whose initial data belong to a bounded set $\Omega$.
We observe that this theorem extends related results in the panorama of literature on the topic as, e.g., \cite[Theorem 2.2 and Remark 2]{CM} where the linear case is considered. 
Further on, as a consequence, we deduce that if the problem has the zero solution, then it is an attractor for all the solutions that originate from a bounded set containing the zero of $E$. 

In Section \ref{s:asm}, in light of the results obtained in abstract spaces, we can resume the study of the robotic model 
subject to a distributed delay. We first present an existence and uniqueness theorem of solutions. 
Afterwards, we reach the goal
to provide conditions for which there is uniform asymptotic stability of the solutions of the equation \eqref{e:pde1} with the boundary conditions \eqref{e:pde2} -\eqref{e:pde3}.

Sections \ref{s:d} and \ref{s:concl} are devoted respectively to the preliminary notions and to the conclusions of the manuscript. In particular, il Section \ref{s:concl} we summarize the main results of the paper and suggest some prospects for future developments of this research.

\section{Preliminaries}
\label{s:d}

Let $E$ be a real Banach space endowed with the norm $\|\cdot\|$ and $J$ a compact interval in $\erre$.
\\
As usual, by $C(J,E)$ we denote the space of $E$-valued continuous functions defined on $J$ endowed with the sup-norm, and by $L^p(J,E)$ the space of all functions $v:J\to E$ such that $v^p$ is Bochner integrable endowed with the norm $\|v\|_{L^p(J,E)}= \left(\int_{J}\|v(z)\|^p\, dz\right)^{\frac{1}{p}}$ (if $E=\erre$, we simply write $L^p(J)$ and $\|v\|_{L^p}$ respectively), $p\ge 1$.
\\
Moreover, by the symbol $L^1_{loc}([a,+\infty[,E)$, $a\in \erre$, we denote the space of all functions $v:[a,+\infty[\to E$ such that $v\in L^1(I,E)$ for every compact interval $I\subset [a,+\infty[$ (shortly, $L^1_{loc}([a,+\infty[)$ if $E=\erre$).
\\
Also, we put $H^k(0,1):=\{y:[0,1]\to \erre: y,y^{(1)},...,y^{(k)}\in L^2([0,1])\}$.
\sm

Let ${\cal L}(E)$ be the Banach space of all bounded linear operators from $E$ to $E$ with the operator norm
$$
\|L\|_{{\cal L}(E)}=\sup_{\|x\|\le 1}\|L x\|
$$
for each $L\in {\cal L}(E)$.
We recall that a family  $\{U(t)\}_{t\ge 0}$ in ${\cal L}(E)$ is said to be a {\em $C_0$-semigroup} (see, e.g. \cite{V}) if
\begin{description}
\item{(U1)} $U(0)=I$,
\item{(U2)} $U(t+s)=U(t)U(s)$, for every $t,s\ge 0$,
\item{(U3)} $\lim_{t\downarrow 0}U(t)x=x$,  for every $x\in E$.
\end{description}
For every $C_0$-semigroup $\{U(t)\}_{t\ge 0}$ it holds that  (see, e.g. \cite[Theorem 2.3.1]{V})
\begin{description}
\item{(U4)} there are constants $D\ge 1, \delta\in \erre$ such that for every $t\ge 0$ we have
\begin{equation}\label{e:D}
\|U(t)\|_{{\cal L}(E)}\le D e^{\delta t}.
\end{equation}
\end{description}
Sometimes in order to put in evidence the involved constants, the $C_0$-semigroup is called {\em of type $(D,\delta)$}. In particular,
\begin{itemize}
\item
if $\{U(t)\}_{t\ge 0}$ is of type $(1,0)$, then it is said to be a {\em $C_0$-semigroup of contractions}; 
\item
if $\{U(t)\}_{t\ge 0}$ is of type $(D,-\omega)$, with $\omega>0$, then it is said to be {\em exponentially stable}.
\end{itemize}

\np
The {\em infinitesimal generator} of a $C_0$-semigroup is the linear operator $A:D(A)\subset E\to E$ defined by
$$
Ax=\lim_{t\downarrow 0}\frac1t (U(t)x-x),
$$
for $x\in D(A):=\left\{ x\in E : \exists \lim_{t\downarrow 0}\frac1t (U(t)x-x) \right\}.$
The set $D(A)$ is dense in $E$ and $A$ is a closed operator (see, e.g. \cite[Theorem 2.4.1]{V}).
\sm

Finally, we recall the Kuratowski and Hausdorff measures of noncompactness on the family of nonempty bounded subsets of $E$, $\alpha$ and $\chi$ respectively, defined as
\begin{eqnarray*}
&& \alpha(\Omega)=\inf\{\varepsilon>0: \Omega \mbox{ can be covered by finitely many sets with diameter } \le \varepsilon\};
\\
&&\chi(\Omega)=\inf\{\varepsilon>0: \Omega \mbox{ can be covered by finitely many balls with radius } \varepsilon\}.
\end{eqnarray*}
These measures of noncompactness are equivalent, since
$\chi (\Omega)\le \alpha (\Omega)\le 2\chi(\Omega)$.
For other properties of the measures of noncompactness we refer, e.g., to \cite{A} or \cite{AT}.
\bg

\section{From the flexible robotics model to the abstract problem}
\label{s:abs}

We consider the fourth-order partial differential equation \eqref{e:pde1} with boundary conditions \eqref{e:pde2}, \eqref{e:pde3} subject to the initial conditions \eqref{e:iv0}-\eqref{e:iv2}.

\np
By using the auxiliary function $\eta:\erre^+_0\to \erre$ defined by (see \cite{CM})
\begin{equation}\label{e:eta}
\eta(t):=-u_{xxx}(t,1)+\dfrac{m}{\beta} u_t(t,1),\ t\ge 0,
\end{equation}
the boundary condition \eqref{e:pde3} can be rewritten as
\begin{equation}\label{e:pde4}
 \eta'(t)+\dfrac1{\beta}\eta(t)+\dfrac1{\beta}\left(\alpha -\dfrac{m}{\beta}\right) u_t(t,1)=0,\ t\ge 0.
\end{equation}
Indeed, by \eqref{e:pde3} and \eqref{e:eta}, we have
$$
mu_{tt}(t,1)+\eta(t)-\dfrac{m}{\beta} u_t(t,1)=-\alpha u_t(t,1)+\beta u_{xxxt}(t,1).
$$
Further, by  
$$
 \beta\eta'(t)=-\beta u_{xxxt}(t,1)+m u_{tt}(t,1),\ t\ge 0,
$$
condition \eqref{e:pde4} holds.
\\
Thus problem  \eqref{e:pde1}-\eqref{e:iv2} is equivalent to the first order system 
\begin{equation*}
\begin{cases}
u_t(t,x)=v(t,x)\\
v_t(t,x)=-u_{xxxx}(t,x)+g\left(t,u(t,x),\int_{0}^t \frac{e^{-(t-s)/T}}{T}\, u(s,x)\, ds\right)\\
\eta'(t)=-\dfrac1{\beta}\eta(t)-\dfrac1{\beta}\left(\alpha -\dfrac{m}{\beta}\right) v(t,1),
\end{cases}
\end{equation*}
with initial conditions
\begin{equation*}
\left\{\begin{array}{l}
u(0,x)=\bar p(x),\\
v(0,x)=\bar q(x),\\
\eta(0)=\bar \eta,
\end{array}\right.
\end{equation*}
and boundary conditions
\begin{equation*}
u(t,0)=u_{x}(t,0)=0,\quad u_{xx}(t,1)=0.
\end{equation*}
This problem can be written as an abstract Cauchy problem in a suitable function space ${\cal H}$,
\begin{equation}
\label{e:H}
{\cal H}:=\{\y :=(p,\ q,\  \eta)^T  : p\in {\cal V},\ q\in L^2([0,1]),\  \eta\in\erre\},
\end{equation}
where  the superscript $T$ stands for the transpose and 
$$
{\cal V}:=\{p\in H^2(0,1): p(0)=p_x(0)=0\},
$$
introducing the linear operator $A:D(A)\subset {\cal H}\to {\cal H}$ defined by
\begin{eqnarray}
\nonumber
&D(A)=\left\{  \y=(p, q,  \eta)^T  : p\in H^4(0,1)\cap {\cal V},\, q\in {\cal V},\  p_{xx}(1)=0,\  \eta=-p_{xxx}(1)+\frac{m}{\beta}q(1)\right\},&
\\
\label{e:A}\\
\nonumber &A\y=A\left(\begin{matrix}p\\ q\\  \eta\end{matrix}\right):=\left(\begin{matrix}q\\ -p_{xxxx} \\  -\frac1{\beta} \eta-\frac1{\beta}\left(\alpha-\frac{m}{\beta}\right)q(1)\end{matrix}\right),&
\end{eqnarray}
and the nonlinear function $f:\erre^+_0\times {\cal H}\times {\cal H}\to {\cal H}$, given by 
\begin{equation}
\label{e:fg}
f(t,\y_1,\y_2)(x):=\left(0_{\cal V}, g(t, p_1(x), p_2(x))  , 0\right)^T,\ x\in [0,1],
\end{equation}
with $g$ sufficiently regular (see Section \ref{s:asm}, condition (g1)).
\\
It is known that ${\cal H}$ is a Hilbert space with a suitable inner product 
 and  that $A$ generates a $C_0$-semigroup on ${\cal H}$, which is at the same time of contractions (see \cite[Theorem 2.1]{CM}) and exponentially stable  (see \cite[Theorem 2.2]{CM}). 
\sm

Therefore, defining $y:\erre_0^+\to  {\cal H}$ as
\beq
\label{e:y}
y(t)(x):=\left(\begin{matrix}u(t,x)\\ v(t,x)\\  \eta(t)\end{matrix}\right), \mbox{\  $t\ge 0,\, x\in [0,1]$,}
\eeq
the problem \eqref{e:pde1}-\eqref{e:iv2}  can be written as the abstract Cauchy problem 
\begin{equation}
\label{e:P_interm}
\begin{cases}
y'(t)= A y(t)+f\left(t,y(t),\int_{t_0}^t  \frac{e^{-(t-s)/T}}{T}y(s)ds\right) ,\ t \ge 0,
\\
y(0)=\bar y,
\end{cases}
\end{equation}
where $\bar y=(\bar p, \bar q, \bar \eta)^T.$
\sm

In the next sections we will deal with this equation in the more general setting given by Banach spaces and with a continuous kernel inside the integral component of the nonlinearity.

\section{Existence, uniqueness and continuous dependence}
\label{s:eudc}

Let $E$ be a real Banach space endowed with a norm $\|\cdot\|$ and $\bv\in E$.  We consider the initial value problem driven by the semilinear integro-differential equation
\begin{equation*}
(P)_{\bv}\ \begin{cases}
y'(t)= A y(t)+f\left(t,y(t),\int_{t_0}^t k(t,s)y(s)ds\right) ,\ t \ge t_0,
\\
y(t_0)=\bv,
\end{cases}
\end{equation*}
where $A:D(A)\to E$ is the infinitesimal generator of a $C_0$-semigroup $\{U(t)\}_{t\ge 0}$ with $\overline{D(A)}=E$, and $f:[t_0,+\infty[\, \times E\times E\to E$, $k:\Delta_\infty:= \{(t,s) \in \erre^2 : t\ge s\ge t_0\} \rightarrow \mathbb{R}^+$
are given functions.\md

A continuous function $y:[t_0,+\infty[ \to E$ is said to be a {\em mild solution}  to $(P)_{\bv}$ if it satisfies the Duhamel formula, i.e.
\begin{equation*}
y(t) = U(t-t_0)\bv+\int_{t_0}^t U(t-s) f\mbox{$\left(s,y(s),\int_{t_0}^s k(s,r)y(r) dr\right)$}\, ds  ,\ t \ge t_0.
\end{equation*}

\np 
The existence of the mild solutions of $(P)_{\bv}$ is guaranteed by the next result, which is a direct consequence of \cite[Corollary 3.1] {R22} and \cite[Remarks 1 and 4]{CR05}.

\begin{proposition}\label{l:exist}
Suppose that the kernel $k$ and the function $f$ satisfy
\begin{description}
\item{(k1)} $k$ is continuous;
\sm

\item{(f1)} for every $v,w \in E$ the map $f(\cdot,v,w)$ is strongly measurable;
\sm

\item{(f2)}
for a.e. $ t \in [t_0,+\infty[$ the map $f(t,\cdot, \cdot)$ is continuous;
\sm

\item{(f3)} there exists a nonnegative function $\mu\in L^1_{loc}([t_0,+\infty[)$ such that,  for a.e. $t \ge t_0$ and all $v,w\in E$,
\begin{equation*}
\|f(t,v,w)\|\le \mu(t)(1+\|v\|+\|w\|)\ ;
\end{equation*}

\item{(f4)} there exists a nonnegative function $h\in L^1_{loc}([t_0,+\infty[)$ such that
\begin{equation*}
\chi(f(t,\Omega_1,\Omega_2))\le h(t)\left[\chi(\Omega_1)+\chi(\Omega_2)\right]\ ,
\end{equation*}
for a.e. $t \ge t_0$ and every bounded $\Omega _1,\Omega _2\subset E$ (where $\chi$ is the Hausdorff measure of noncompactness).
\end{description}
Then problem $(P)_{\bv}$ has at least one mild solution on $[t_0,+\infty[$.
\end{proposition}

Under stronger assumptions on the nonlinearity $f$ we obtain the uniqueness of the mild solution and the continuous dependence on the initial data.
\begin{theorem}\label{c:e}
Assume that the kernel $k$ has property (k1) and the function $f$ satisfies (f1),
\begin{description}
\item{(f5)}
there exists $C>0$ such that for a.e. $t \ge t_0$
the map $f(t,\cdot, \cdot)$ is $C$-lipschitzian, i.e. 
$$
\|f(t,v_1,v_2)-f(t,w_1,w_2)\|\le C(\|v_1-w_1\|+\|v_2-w_2\|),
$$
for all $v_1,w_1,v_2,w_2\in E$;
\item{(f6)} $\|f(\cdot, 0,0)\|\in L^1_{loc}([t_0,+\infty[)$.
\end{description}
Then, for every $\bv\in E$  the problem $(P)_{\bv}$ has a unique mild solution, continuously depending on the initial data.
\end{theorem}
\np {\bf Proof.}
First of all, hypothesis (f5) trivially implies property (f2). Furthermore, by (f5) and (f6) also property (f3) holds, since  for almost all $t\ge t_0$ and $v,w\in E$ we get
$$
\|f(t,v,w)\|\le \|f(t,0,0)\|+C(\|v\|+\|w\|) \le (\|f(t,0,0)\|+C)(1+\|v\|+\|w\|),
$$
and taking $\mu(\cdot):=\|f(\cdot,0,0)\|+C$ property (f3) is satisfied.
\\
Finally, by (f5) it follows (f4) as well. In fact, for a.e. $t\ge t_0$ such that the map $f(t,\cdot, \cdot)$ is $C$-lipschitzian,  let us consider
\begin{eqnarray*}
&&\alpha(f)(t) := \inf\{k(t)>0 : \alpha(f(t,\Omega_1,\Omega_2))\le k(t)[\alpha(\Omega_1)+\alpha(\Omega_2)] \mbox{ for bounded } \Omega_1,\Omega_2\subset E\}\\
&&Lip(f)(t)    := \inf\{k(t)>0 : \|f(t,v_1,v_2)-f(t,w_1,w_2)\|\le k(t)(\|v_1-w_1\|+\|v_2-w_2\|) \},
\end{eqnarray*}
where $\alpha$ and $\chi$ are the Kuratowski and Hausdorff measures of noncompactness respectively.
For all bounded $ \Omega_1,\Omega_2\subset E$, we get (see \cite[Section 2]{AEFV})
\begin{eqnarray*}
\chi(f(t,\Omega_1,\Omega_2))&\le& \alpha(f(t,\Omega_1,\Omega_2)) 
\le \alpha(f)(t)[\alpha(\Omega_1)+\alpha(\Omega_2)]\\
&\le& Lip(f)(t)  [\alpha(\Omega_1)+\alpha(\Omega_2)] \le  2C  [\chi(\Omega_1)+\chi(\Omega_2)],
\end{eqnarray*}
and (f4) follows just taking $h(t)=2C$, a.e. $t\ge t_0$.
\\
Hence we can apply Proposition \ref{l:exist} and claim that for every fixed $\bv\in E$  the solution set of  problem  $(P)_{\bv}$ is nonempty. 
\\
Let us show that in the present setting $(P)_{\bv}$ is a singleton. Let $y,z\in C([t_0,+\infty[,E)$ be two mild solutions of $(P)_{\bv}$. 
Let us fix $n\in \enne$ with $n>t_0$. 
Then, for every $t\in [t_0,n]$ we have
\begin{eqnarray*}
\|y(t)-z(t)\|
&\le& \int_{t_0}^t \|U(t-s)\|_{{\cal L}(E)} \left\|f\left(s,y(s),\int_{t_0}^sk(s,r)y(r)dr\right)-f\left(s,z(s),\int_{t_0}^sk(s,r)z(r)dr\right)\right\| ds.
\end{eqnarray*}
Clearly from (U4) (see Section \ref{s:d}), there exists $M_n>0$ such that $\|U(t-s)\|_{{\cal L}(E)}\le M_n$ for every $(t,s)\in \Delta_n:=\{(t,s)\in\erre^2:t_0\le s\le t\le n\}$. Hence, by using  (f5), we get
\begin{eqnarray*}
\|y(t)-z(t)\|
&\le& M_nC\int_{t_0}^t \|y(s)-z(s)\|ds+M_nC\int_{t_0}^t \int_{t_0}^s k(s,r)\|y(r) - z(r)\|dr ds.
\end{eqnarray*}
By the continuity of the function $k$ there exists $K_n>0$ such that $k(t,s)\le K_n$ for every $(t,s)\in \Delta_n$. Therefore, by changing the order of the integrals and renaming the variables, we deduce
\begin{eqnarray*}
\|y(t)-z(t)\|
&\le& M_nC\int_{t_0}^t \|y(s)-z(s)\|ds+M_nCK_n\int_{t_0}^t\int_{r}^t \|y(r) - z(r)\|ds dr\\
&\le& M_nC\int_{t_0}^t \|y(s)-z(s)\|ds+M_nCK_n(n-t_0) \int_{t_0}^t\|y(r) - z(r)\|dr\\
&=& M_nC[1+K_n(n-t_0) ]\int_{t_0}^t \|y(s)-z(s)\|ds.
\end{eqnarray*}
The Gronwall inequality yields $\|y(t)-z(t)\|=0$ for every $t\in [t_0,n]$, and from the arbitrariness of $n$ it follows that  $\|y(t)-z(t)\|=0$ for every $t\ge t_0$.

\medskip
Finally, let us prove the continuous dependence of the mild solutions on the initial data. To this aim we define  
the function $\varphi:[t_0,+\infty[ \times E\to E$ as
$$
\varphi (t,\bv):=y(t) , \ t\ge t_0,
$$
where $y$ is the unique mild solution to $(P)_{\bv}$.
\\
Let $(t,\bv)\in [t_0,+\infty[ \times E$ be arbitrarily fixed. Of course, there exists $n>t_0$ such that $t\in [t_0,n[$. Now, we consider any $(s,\bw)\in [t_0,n] \times E$. Denoted by $y,z$ the unique mild solutions respectively to $(P)_{\bv}$ and  $(P)_{\bw}$, we have
\begin{eqnarray*}
\|\varphi (t,\bv)-\varphi (s,\bw)\| & =&\|y(t)-z(s)\|\,\, \le \,\, 
\|y(t)-z(t)\|+\|z(t)-z(s)\|.
\end{eqnarray*}
With the same arguments as above, we achieve the following estimate
\begin{eqnarray*}
\|\varphi (t,\bv)-\varphi (s,\bw)\| 
&\le&M_n\|\bv-\bw\|+M_nC[1+K_n(n-t_0)]\int_{t_0}^t\|y(s)-z(s)\|ds\\
&&+\|z(t)-z(s)\|.
\end{eqnarray*}
By using again the Gronwall inequality we obtain
\begin{eqnarray*}
\|\varphi (t,\bv)-\varphi (s,\bw)\| 
&\le&M_n\|\bv-\bw\|e^{M_nC[1+K_n(n-t_0)](n-t_0)}+\|z(t)-z(s)\|.
\end{eqnarray*}
Therefore, passing to the limit for $(s,\bw)\to (t,\bv)$ we have $\|\varphi (s,\bw)-\varphi (t,\bv)\|\rightarrow 0$, i.e. $\varphi $ is continuous in $(t,\bv)$. The arbitrariness of $(t,\bv)$ concludes the proof.
\hfill$\Box$
\bg

Generalizing some classical results (see for instance \cite[Theorem 1.6]{P}), if $f$ satisfies a Lipschitz condition in the three variables, the kernel $k$ is Lipschitz in the first variable and the initial datum belongs to the domain of the operator $A$, then the corresponding mild solution of $(P)_{\bv}$ is a \emph{strong} solution, i.e, it is differentiable almost everywhere on $[t_0, +\infty[$ with derivative in $L^1_{loc}([t_0,+\infty[, E)$, 
and satisfies the initial condition and the differential equation in $(P)_{\bv}$ almost everywhere on $[t_0, +\infty[$. More precisely, the following holds.
\begin{theorem}\label{m>s}
Suppose that the kernel $k$ satisfies (k1) and 
\begin{description}
\item{(k2)} for every $n \in \mathbb N, \, n > t_0$ there exists $\tilde k_n>0$ such that 
$$|k(t_1,s)-k(t_2,s)|\leq \tilde{k}_n |t_1-t_2|,$$
for all $(t_1,s), (t_2,s) \in \Delta_n=\{(t,s)\in\erre^2:t_0\le s\le t\le n\}$.
\end{description}
\noindent Moreover, assume that $f$ satisfies properties (f5) and
\begin{description}
\item{(f7)}
for every $n \in \mathbb N, \, n > t_0$ there exists $C_n>0$ such that 
$$
\|f(t_1,v,w)-f(t_2,v,w)\|\le C_n|t_1-t_2|,
$$
for all  $t_1, t_2 \in [t_0, n]$, \, $v,w\in E$.
 \end{description}
Then, for every $\bv\in D(A)$  the problem $(P)_{\bv}$ has a unique strong solution on $[t_0, \infty[$, continuously depending on the initial data.
\end{theorem}

\np {\bf Proof.} First of all, notice that assumption (f7) trivially implies (f1) and (f6) of Theorem \ref{c:e}. Therefore $(P)_{\bv}$ has a unique mild solution  $y \in C([t_0, \infty[, E)$, continuously depending on the initial data.  Let $\psi:[t_0, \infty[ \to  E $ be  defined by
\[\psi(t)=\int_{t_0}^t k(t,s) y(s)\, ds.\]
If  $n \in \mathbb{N}, n > t_0$ is fixed,  the assumption (k2) implies that the function $\psi$ is Lipschitz-continuous on $[t_0, n]$. Indeed, for every $t, t+h \in [t_0, n]$, we have
\begin{align*}
\| \psi(t+h)-\psi(t)\| &\leq \int_{t_0}^t |k(t+h, s)-k(t,s)| \| y(s) \| \, ds+ \left| \int_t^{t+h} k(t+h, s) \| y(s) \| \, ds \right|\\
&\leq \tilde k_n |h| (n-t_0) \bar y_n+ K_n |h| \bar y_n := \hat K_n |h|,
\end{align*}
where $\bar y_n:= \max_{s\in [t_0, n]} \| y(s) \|$, and $K_n:= \max_{(t,s)\in \Delta_n} k(t,s)$. 
\\
Now we prove that also $y$ is Lipschitz-continuous on $[t_0, n]$. For $t, t+h \in [t_0, n[$, let us consider
 \begin{align*}
 y(t+h)-y(t)= &[U(t+h-t_0)-U(t-t_0)] \bv + \int_{t_0}^{t_0+h} U(t+h-s)f(s, y(s), \psi(s))\, ds\\
 &+\int_{t_0+h}^{t+h} U(t+h-s)f(s, y(s), \psi(s))\, ds- \int_{t_0}^{t} U(t-s)f(s, y(s), \psi(s))\, ds\\
 =& \int_{t-t_0}^{t+h-t_0} U(s) A \bv \, ds + \int_{t_0}^{t_0+h} U(t+h-s)f(s, y(s), \psi(s))\, ds\\
 &+ \int_{t_0}^{t} U(t-s)\left[f(s+h, y(s+h), \psi(s+h)) - f(s, y(s), \psi(s))\right]\, ds,
 \end{align*}
 where we used a basic property of $C_0$-semigroups, see, e.g., \cite[Theorem 2.4, d)]{P}.
 Then, put $M_n= \max_{(t,s)\in \Delta_n} \|U(t-s)\|_{{\cal L}(E)}$, and  $F_n:= \max_{t \in [t_0, n]} f(t, y(t), \psi(t))$, and taking into account the Lipschitz-continuity of $\psi$, and (f5), (f7), we get
 \begin{align*}
 \|y(t+h)&-y(t)\|\leq M_n |h| \|A \bv\|+ M_n |h| F_n\\
 &+M_n  \int_{t_0}^t (C_n|h| + C\|y(s+h)-y(s)\| + C\|\psi(s+h)-\psi(s)\|)\, ds\\
 & \leq  |h| M_n \big(\|A \bv\|+  F_n + C_n (n-t_0)+C\hat K_n(n-t_0) \big)+ M_n C \int_{t_0}^t \|y(s+h)-y(s)\| \, ds.
\end{align*}
 Put $\tilde C_n:=M_n  \big(\|A \bv\|+  F_n + C_n (n-t_0)+C\hat K_n(n-t_0)\big)$, the Gronwall inequality yields
 \[
  \|y(t+h)-y(t)\|\leq |h| \tilde C_n e^{(n-t_0)M_n C},
 \]
i.e.,  $y$ is Lipschitz-continuous on $[t_0, n]$. Thus, by composition, also the function 
\[
\tilde f(t):= f(t, y(t), \psi(t)), \quad t \geq t_0
\]
is Lipschitz-continuous on $[t_0, n]$, and the Cauchy problem 
\[
\begin{cases}
z'(t)= A z(t) + \tilde f(t), &t \in [t_0, n]\\
z(t_0)=\bv
\end{cases}
\]
has a unique strong solution $z_{\bv}$ on $[t_0, n]$ (see \cite[Corollary 4.2.11]{P}), which satisfies
\begin{align*}
z_{\bv}(t)&=U(t-t_0) \bv +\int_{t_0}^t U(t-s) \tilde f(s) \, ds\\
&=U(t-t_0) \bv +\int_{t_0}^t U(t-s) f(s, y(s), \psi(s) )\, ds = y(t)
\end{align*}
for all $t \in [t_0, n]$, i.e., $y$ is a strong solution of $(P)_{\bv}$ on $[t_0, n]$. From the  arbitrariness of $n$ follows the thesis.
\hfill$\Box$

\section{Uniform asymptotic stability and attractors}
\label{s:as}

In the line of \cite{BEA,BR,CM}, we give the following definition of local uniform asymptotic stability of mild solutions to the driving equations of  problems $(P)_{\bv}$ (see also \cite{BD,R21}).

\begin{definition}
Let $\Omega$ be a nonempty bounded subset of $E$. The mild solutions of the integro-differential equation \eqref{e:eq}
are said to be {\em uniformly asymptotically stable}  on  $\Omega$ if
\sm

\begin{itemize}
\item[(as)]
for every $\varepsilon>0$ there exists $t(\varepsilon)>0$ such that
$$
\|z(t)-y(t)\|\le \varepsilon, \mbox{ for every $t\ge t(\varepsilon)$}
$$
for every $z,y\in C([t_0,+\infty[,E)$ mild solutions of $(P)_{\bw}, (P)_{\bv}$ respectively, for every $\bw, \bv\in\Omega$.
\end{itemize}
\end{definition}
In other words, the mild solutions of \eqref{e:eq} are uniformly asymptotically stable on $\Omega$ if the limit
$$
\lim_{t\to +\infty}\|z(t)-y(t)\|=0
$$
is uniform in the set
\begin{equation}\label{e:S}
{\cal S}(\Omega)=\left\{ y: y \mbox{ is a mild solution of } (P)_{\bv},\, \bv\in\Omega \right\}.
\end{equation}


\sm

\np
The next result provides a sufficient condition for the uniform asymptotic stability.

\begin{theorem}\label{t:as1}
Suppose that  the linear operator $A$ generates an  exponentially stable $C_0$-semigroup $\{U(t)\}_{t\ge 0}$ of type $(D,-\omega)$ with $\omega>0$, and that the function $f$ has properties (f1), (f5), (f6).
Assume that the function $k$ satisfies (k1) and
\begin{description}
\item{(k3)} there exist $a,b>0$  such that
$$
k(t,s)\le a e^{-b (t-s)},\ (t,s)\in \Delta_\infty.
$$
\end{description}
Suppose also that the next inequalities hold:
\begin{equation}
\label{(i)} 
C<\frac{\omega}{D}, \quad b>\omega, \quad a<\frac{(b-\omega)(\omega-CD)}{CD}. 
\end{equation}
\np
Then, for every bounded $\Omega\subset E$  the mild solutions of equation \eqref{e:eq} are  uniformly asymptotically stable on $\Omega$.
\end{theorem}

\np {\bf Proof.}
By Theorem \ref{c:e},  for every $\bv \in E$ the problem  $(P)_{\bv}$ has a unique mild solution.
So, for $\bw, \bv\in\Omega$ arbitrarily fixed, let  $z$ and $y$ be the mild solutions of $(P)_{\bw}$ and $(P)_{\bv}$ respectively.

\np
For every $t\ge t_0$, by (U4),  (f5), we obtain
\begin{eqnarray*}
\|z(t)-y(t)\|&\le &
\left\|U(t-t_0)(\bw-\bv)\right\| \\
&&+\int_{t_0}^t \left\|U(t-s)\left[f\mbox{$\left(s,z(s),\int_{t_0}^s k(s,r)z(r) dr\right)$} - f\mbox{$\left(s,y(s),\int_{t_0}^s k(s,r)y(r) dr\right)$}\right]\right\|\, ds\\
&\le &  De^{-\omega(t-t_0)}\left\|\bw-\bv\right\| \\
&&+\int_{t_0}^t De^{-\omega(t-s)}\left\|f\mbox{$\left(s,z(s),\int_{t_0}^s k(s,r)z(r) dr\right)$} - f\mbox{$\left(s,y(s),\int_{t_0}^s k(s,r)y(r) dr\right)$}\right\|\, ds  \\
&\le &  De^{-\omega(t-t_0)}\left\|\bw- \bv\right\| \\
&&+\int_{t_0}^t De^{-\omega(t-s)}C \left[ \|z(s)-y(s)\|
+\left\| \int_{t_0}^s k(s,r)z(r) dr - \int_{t_0}^s k(s,r)y(r) dr \right\| \right] \, ds\\
&\le &  De^{-\omega(t-t_0)}\left\|\bw- \bv\right\| \\
&&+\int_{t_0}^t C De^{-\omega(t-s)}\left[ \|z(s)-y(s)\|
+\int_{t_0}^s k(s,r)\left\| z(r) -y(r)  \right\|dr \right] \, ds.
\end{eqnarray*}
Therefore,
\begin{eqnarray*}
\|z(t)-y(t)\|e^{\omega(t-t_0)}&\le &  D\left\|\bw- \bv\right\| \\
&&+C D\int_{t_0}^t  \|z(s)-y(s)\|e^{\omega(s-t_0)}\, ds\\
&&+C D\int_{t_0}^t\left( \int_{t_0}^s k(s,r)\left\| z(r) -y(r)  \right\|e^{\omega(s-t_0)}dr \right)\, ds
\end{eqnarray*}
Put
\begin{equation}
\label{e:gamma}
\gamma(t):=\|z(t)-y(t)\|e^{\omega(t-t_0)} \mbox{ and }
c_0:=D\left\|\bw- \bv\right\|,
\end{equation}
the previous inequality reads as
\begin{eqnarray*}
\gamma(t) &\le &  c_0
 +C D\int_{t_0}^t  \gamma(s)\, ds
+C D\int_{t_0}^t \left(\int_{t_0}^s k(s,r) \gamma(r) e^{\omega(s-r)} dr\right) \, ds.
\end{eqnarray*}
Since, by changing the order of the integrals and renaming the variables, it holds that
\begin{eqnarray*}
   \int_{t_0}^t \left(\int_{t_0}^s k(s,r) \gamma(r) e^{\omega(s-r)} dr\right) \, ds
&=& \int_{t_0}^t\left( \int_{r}^t k(s,r)  e^{\omega(s-r)} ds \,\right) \gamma(r)\, dr\\
&=& \int_{t_0}^t  \left(\int_{s}^t k(\sigma,s)  e^{\omega(\sigma-s)} d\sigma\,\right) \gamma(s)\, ds,
\end{eqnarray*}
then we get
\begin{equation}
\label{e:G}
\gamma(t) \le
c_0
+\int_{t_0}^t  C D\gamma(s)\, ds
+ \int_{t_0}^t C D  \left(\int_{s}^t k(\sigma,s)  e^{\omega(\sigma-s)} d\sigma\,\right) \gamma(s) \, ds.
\end{equation}
Now, by  (k3) we deduce that
\begin{eqnarray}
\nonumber \int_{s}^t k(\sigma,s)  e^{\omega(\sigma-s)} d\sigma&\le&
\int_{s}^t a e^{(\omega-b)(\sigma-s)} d\sigma\\
\label{e:k1-cons} &=&
-\dfrac{a}{\omega-b}\left( 1-e^{(\omega-b)(t-s)} \right) <\dfrac{a}{b-\omega}.
\end{eqnarray}
Hence, by \eqref{e:G} and \eqref{e:k1-cons} we have
$$
 \gamma(t) < c_0
+ \int_{t_0}^t  C D \left(1+\dfrac{a}{b-\omega}\right) \gamma(s) \, ds.
$$
We can therefore apply the Gronwall inequality,  obtaining for all $t\ge t_0$
$$
\gamma(t) \le c_0 \, e^{C D\left(1+\frac{a}{b-\omega}\right) (t-t_0)}.
$$
By the definitions of $\gamma$ and $c_0$ (see \eqref{e:gamma}) we get
\begin{equation}
\label{e:stima_finale}
\|z(t)-y(t)\|
\le D \|\bw-\bv\|
 \, e^{\left[C D\left(1+\frac{a}{b-\omega}\right) -\omega\right](t-t_0)},\, t\ge t_0.
\end{equation}
Since hypothesis \eqref{(i)} implies that
$$
C D\left(1+\frac{a}{b-\omega}\right) -\omega<0,
$$
from \eqref{e:stima_finale} we have
\begin{equation}
\label{e:lim}
\lim_{t\to +\infty}\|z(t)-y(t)\|=0.
\end{equation}
The boundedness of $\Omega$ gives the uniformity of the limit in ${\cal S}(\Omega)$, i.e. the uniform asymptotic stability on $\Omega$ of the mild solutions of \eqref{e:eq}.
\hfill $\Box$

\begin{remark}
Notice that, from the proof of Theorem \ref{t:as1}, it follows that the norm of the difference of two mild solutions of equation \eqref{e:eq} decays to $0$ exponentially, even if $\Omega$ is not bounded.
Thus our theorem extends some related results in the literature (see for instance \cite[Theorem 2.2 and Remark 2]{CM}, where the case of the linear equation \eqref{e:pde0} is considered).  
\end{remark} 

\np
Following \cite{DL}, we set the next definition.
\begin{definition}
A line $y(t)=c$, where $c\in E$, is called an {\em attractor} for a mild solution $z$ of the integro-differential equation \eqref{e:eq} if \
\begin{equation*}
\lim_{t\to+\infty}\|z(t)-c\|=0.
\end{equation*}
\end{definition}
By strenghtening  hypothesis (f6), we achieve the following result on the attractors and on the set ${\cal S}(\Omega)$ (cf. \eqref{e:S}).


\begin{corollary}\label{c:attr}
Let $\Omega\subset E$ a bounded set with $0\in \Omega$.
Suppose that  the hypotheses of Theorem \ref{t:as1} are satisfied and that the function $f$ has property
\begin{description}
\item{(f6*)} $f(\cdot, 0,0)=0$.
\end{description}

\np
Then, the solution $y_0\equiv 0$ of \eqref{e:eq} is an attractor for all the mild solutions belonging to ${\cal S}(\Omega)$, and ${\cal S}(\Omega)$ is bounded.
\end{corollary}
\np {\bf Proof.}
Clearly the null function $y_0$ is a strong solution of \eqref{e:eq} by (f6*).
Moreover, $y_0$ is also a mild solution of $(P)_0$. Since $0\in \Omega$, if $z_{\bv}$ is the unique mild solution of $(P)_{\bv}$,  $\bv\in \Omega$, retracing the proof of Theorem \ref{t:as1} we obtain (see \eqref{e:lim})
$$
\lim_{t\to +\infty} \|z_{\small \bv}(t)-y_0(t)\|=0,
$$
i.e. the line $y_0$ is an attractor for all the mild solutions belonging to ${\cal S}(\Omega)$.
\\
Finally, by \eqref{e:stima_finale} we can write
$$
\|z_{\bv}(t)\|\le D\,  \mbox{diam} (\Omega),
$$
for every $t\in [t_0,+\infty[$ and $\bv\in\Omega$, that is ${\cal S}(\Omega)$ is bounded.
\hfill $\Box$

\begin{remark}
\label{r:strong-as}
We point out that if we add hypotheses (k2) and (f7) (see Theorem \ref{m>s})  in Theorem \ref{t:as1} and Corollary \ref{c:attr}, then their statements hold for strong solutions of \eqref{e:eq}. 
\end{remark}

\section{Uniform asymptotic stability of the flexible robotics model}
\label{s:asm}

In this section we apply the abstract results developed in Sections \ref{s:eudc} and \ref{s:as}
to the flexible robotics model  \eqref{e:pde1}-\eqref{e:iv2}.

 We first provide the existence of a unique solution to the model as an application of Theorem \ref{m>s}.

\begin{theorem}
\label{l:g-f}
Assume that the nonlinearity $g$ satisfies properties
\begin{description}
\item{(g1)}
for every $t\ge 0$, the map $x\mapsto g(t, p_1(x), p_2(x))$ belongs to $L^2([0,1])$, for every $p_1,p_2\in {\cal V}$;

\item{(g2)}
there exists $C>0$ such that,  for every $t \ge 0$
$$
\|g(t, p_1(\cdot), p_2(\cdot))-g(t, \hat p_1(\cdot), \hat p_2(\cdot))\|_{L^2}
\le C\left( \|p_1-\hat p_1\|_{L^2} + \|p_2-\hat p_2\|_{L^2} \right),
$$
 for all
$p_1,\hat p_1, p_2, \hat p_2\in {\cal V}$;

\item{(g3)} for every $n \in \mathbb N^+$ there exists $C_n>0$ such that 
$$\|g(r,p_1(\cdot), p_2(\cdot))-g(s,p_1(\cdot), p_2(\cdot))\|_{L^2} \le C_n|r-s|,$$
for all $r,s\in [0,n]$, \ $p_1,p_2\in {\cal V}$.
\end{description}
Then, for every $(\bar p,\bar q,\bar\eta)^T\in D(A)$  there exists a unique strong solution $u:\erre_0^+\, \times [0,1]\to \erre$ of   \eqref{e:pde1}-\eqref{e:iv2}  such that
 $u(\cdot,x)$ is continuously differentiable, $u_t(\cdot, x)$ differentiable a.e. on $\erre^+_0$, and $u_{tt}(\cdot, x)\in L^1_{loc}(\erre^+_0)$ 
for every $x\in [0,1]$; $u(t,\cdot)\in H^4(0,1)$,  $u_t(t,\cdot)\in H^2([0,1])$.
\end{theorem}

\np {\bf Proof.} Let us consider the problem \eqref{e:P_interm} in the space $\cal H$ (see \eqref{e:H}). 
We immediately note that the function $f$ defined in \eqref{e:fg} is well-posed by hypothesis (g1).
Moreover, by (g2) for any $\y_1,\y_2,\hat\y_1,\hat\y_2\in {\cal H}$ we get
\begin{eqnarray*}
\|f(t,\y_1,\y_2)-f(t,\hat\y_1,\hat\y_2)\|_{\cal H} &= &
\|g(t, p_1(\cdot), p_2(\cdot))-g(t, \hat p_1(\cdot),\hat p_2(\cdot))\|_{L^2}\\
&\le &  C\left( \|p_1-\hat p_1\|_{L^2} + \|p_2-\hat p_2\|_{L^2} \right)\\
&\le &  C\left( \|\y_1-\hat\y_1\|_{{\cal H}} +\|\y_2-\hat\y_2\|_{{\cal H}} \right),
\end{eqnarray*}
so that (f5) holds for $f$.
Similarly, (g3) implies (f7).
\\
Since in the model the kernel is $k(t,s)=\frac{e^{-(t-s)/T}}{T}$, it is easily seen that it satisfies both (k1) and (k2) with $\tilde k_n=\frac1{T^2}$ for all $n\in \enne^+$. Then we can apply Theorem \ref{m>s} and claim that \eqref{e:P_interm} has a unique strong solution
$y:\erre_0^+\, \to D(A)\subset {\cal H}$ such that
$$
y'(t)=Ay(t)+ f\mbox{$\left(t,y(t),\int_{0}^t \frac{e^{-(t-s)/T}}T y(s) ds\right)$}, \ a.e.\ t\ge 0
$$
and \[y(0)=(\bar p,\bar q,\bar\eta)^T.\]
According to \eqref{e:y}  we achieve the existence of a unique  solution $u:\erre_0^+\, \times [0,1]\to \erre$ of the  flexible robotics model  satisfying the thesis.
\hfill $\Box$
\sm

\begin{remark}
Clearly  Theorem \ref{m>s} guarrantees more, that is also the continuous dependence of the solutions of the model by the initial data. Further, notice that the auxiliary function $\eta$ defined in \eqref{e:eta} is a $C^1$-function because it satisfies \eqref{e:pde4}.
\end{remark}

We can now state the uniform asymptotic stability result for the solutions of equation \eqref{e:pde1} subject to the boundary conditions \eqref{e:pde2}-\eqref{e:pde3}.

\begin{theorem}
\label{t:appl}
Assume that the nonlinearity $g$ satisfies properties (g1)-(g3).
\\
If the constants $C$ and $T$, which appear respectively in (g2) and in the kernel $k$, are sufficiently small, then
for every bounded $\Omega\subset D(A)
$ the strong solutions of  \eqref{e:pde1}-\eqref{e:iv2} with $(\bar p, \bar q,\bar \eta)^T\in\Omega$ are uniformly asymptotically stable.
\end{theorem}

\np {\bf Proof.}
As already observed in Section \ref{s:abs}, the linear operator $A$ defined in \eqref{e:A}
 generates a $C_0$-semigroup on ${\cal H}$ which is both of contractions and exponentially stable. So the property \eqref{e:D} holds for $D=1$ and $\delta=-\omega$, for some $\omega>0$.
Moreover, as in the previous theorem, from (g1)-(g3) properties (f5), (f7) follow. 
Being the kernel $k(t,s)=\frac{e^{-(t-s)/T}}{T}$, then it satisfies (k1), (k2) and even property (k3), taking $a=b=\frac1T$.
Further, a simple computation shows that \eqref{(i)} holds, just taking 
\begin{equation}
\label{e:CT}
C<\frac{\omega}2 \quad \mbox{ and } \quad  T<\frac{\omega-2C}{\omega(\omega-C)}. 
\end{equation}  
\\
Therefore all the assumptions of Theorem \ref{t:as1} are satisfied if $C$ and $T$ are sufficiently small. Hence, taking into account  Remark \ref{r:strong-as}, we can conclude that for every bounded set $\Omega\subset D(A)$
the solutions of system \eqref{e:pde1}-\eqref{e:pde3} with initial data in $\Omega$ are uniformly asymptotically stable.
\hfill $\Box$
\bg

\begin{example}
Let us consider the equation
\begin{equation}
\label{e:ex}
 u_{tt}(t,x)+\, u_{xxxx}(t,x) + \gamma^2 u(t,x)+\lambda \int_{0}^t \frac{e^{-(t-s)/T}}{T}\, u(s,x)\, ds=0,\ t\ge 0, \, x\in [0,1],
\end{equation}
with the boundary conditions  \eqref{e:pde2}, \eqref{e:pde3}, where $\gamma \in \erre$, $\lambda \in \erre^+_0$. Assume that $C:=\max\{\lambda,\gamma^2\}$ and $T$ satisfy \eqref{e:CT}.
Then, according to Theorem \ref{t:appl}, for every bounded set $\Omega\subset D(A)$ the strong solutions of problem  \eqref{e:ex}, \eqref{e:pde2}, \eqref{e:pde3}  are uniformly asymptotically stable on $\Omega$.
\sm

\np
Indeed, the above equation is obviously a particular case of \eqref{e:pde1}, just taking
$$
g\left(t,u(t,x),\int_{0}^t \frac{e^{-(t-s)/T}}{T}\, u(s,x)\, ds\right):=-\gamma^2 u(t,x)-\lambda \int_{0}^t \frac{e^{-(t-s)/T}}{T}\, u(s,x)\, ds,
$$
and it is easy to check that the function
$$
g(t,p_1,p_2):=-\gamma^2 p_1-\lambda p_2,\, t\in \erre_0^+,\, p_1,p_2\in \erre,
$$
satisfies properties (g1)-(g4).
\sm

Notice that in the case $\lambda =0$ equation \eqref{e:ex} becomes the equation of the linear beam equation 
$$
 u_{tt}(t,x)+\, u_{xxxx}(t,x) + \gamma^2 u(t,x)=0,\ t\ge 0, \, x\in [0,1],
$$
studied, for instance, in \cite{GG} with different boundary conditions.

\end{example}

\section{Conclusions}
\label{s:concl}

In this article we have shown the existence, uniqueness and continuous dependence on initial data for solutions of a boundary value problem for a system od higher order PDEs arising from a flexible robotics model. Further, the uniform asymptotic stability of solutions and the existence of attactors have been proved if the nonlinear term has small displacements and the fading delay has a small width. The results were achieved thanks to the transformation of the problem into an abstract Cauchy problem for a semilinear integro-differential equation in Banach spaces.
\\
We believe that this article lends itself to possible further investigations.
First of all, we expect that the abstract results of Sections \ref{s:eudc} and \ref{s:as} can also be applied to other models, such as those of flexible beams studied in \cite{C} or \cite{GG}. Furthermore, from the purely mathematical analysis point of view, we believe that it is possible to deepen the study of uniform asymptotic stability for integro-differential equations by taking into consideration, for example, problems referable to \eqref{e:P_interm}, even in presence of external impulsive forces. 
\\
Numerical simulations could be done in a future work in order to estimate upper bounds for the constants $C$, $T$ (see Theorem \ref{t:appl}) and to illustrate the theoretical findings.
\bg

\np
{\bf Acknowledgement and funding.} 

This study was partly funded by: Research project of MIUR (Italian Ministry of Education, University and Research) Prin 2022  “Nonlinear differential problems with applications to real phenomena” (Grant Number: 2022ZXZTN2). 

The authors are members of the national group GNAMPA of INdAM 
and, one of them, of the UMI Group T.A.A. “Approximation Theory and Applications”.
\md

\np
{\bf Conflicts of interest.} The authors declare no conflict of interest.
\md

\end{document}